\documentclass[letter,12pt]{article}

\newlength{\xtextwidth}
\newlength{\xtextheight}
\newlength{\xmargin}
\usepackage[text={\xtextwidth,\xtextheight},centering]{geometry}

\usepackage{lscape}

\usepackage{graphicx}
\usepackage{graphics}
\usepackage{epsfig}
\usepackage{amsfonts}
\usepackage{amsmath}
\usepackage{amssymb}
\usepackage{amsthm}
\usepackage{fancyhdr}
\usepackage{indentfirst}
\usepackage{algorithm}
\usepackage{algorithmic}
\usepackage{multirow}
\usepackage{xcolor}
\usepackage{color}
\usepackage{listings}
\everymath{\displaystyle}

\newcommand{\closure}[2][3]{{}\mkern#1mu\overline{\mkern-#1mu#2}} 

\newtheorem{example}{Example}
\theoremstyle{definition} 
\theoremstyle{plain} 
\theoremstyle{plain} 
\theoremstyle{remark} 
\theoremstyle{plain} \swapnumbers

\numberwithin{algorithm}{section}

\begin{document}
\title{A Parallel Auxiliary Grid AMG Method for GPU}

\date{}

\author{Lu Wang, \thanks{Department of Mathematics, The Pennsylvania State University, University Park, PA 16802. Email:
    \texttt{wang\_l@math.psu.edu}}
\and Xiaozhe Hu,  \thanks{Department of Mathematics, The Pennsylvania State University, University Park, PA 16802. Email:
    \texttt{hu\_x@math.psu.edu}}
    \and Jonathan Cohen, \thanks{NVIDIA Research, Email: 
    \texttt{jocohen@nvidia.com}}
   \and Jinchao Xu  \thanks{Department of Mathematics, The Pennsylvania State University, University Park, PA 16802. Email:
    \texttt{xu@math.psu.edu}}
 }
 
\maketitle

\begin{abstract}
In this paper, we develop a new parallel auxiliary grid algebraic multigrid (AMG) method to leverage the power of graphic processing units (GPUs).  In the construction of the hierarchical coarse grid, we use a simple and fixed coarsening procedure based on a region quadtree generated from an auxiliary grid. This allows us to explicitly control the sparsity patterns and operator complexities of the AMG solver. This feature provides (nearly) optimal load balancing and predictable communication patterns, which makes our new algorithm suitable for parallel computing, especially on GPU. We also design a parallel smoother based on the special coloring of the quadtree to accelerate the convergence rate and improve the parallel performance of this solver.  Based on the CUDA toolkit~\cite{Nvidia2010}, we implemented our new parallel auxiliary grid AMG method on GPU and the numerical results of this implementation demonstrate the efficiency of our new method. The results achieve an average speedup of over 4 on quasi-uniform grids and 2 on shape regular grids when compared to the AMG implementation in CUSP \cite{Garland2010}.

\end{abstract}

\vspace{1ex} 


\section{Introduction} 
We consider the iterative solvers for large-scale sparse linear systems 
\begin{equation}
Ax = b.
\end{equation} 
that arise from the discretizations of the partial differential equations (PDEs).  The multigrid (MG) method is one of the most efficient iterative solvers. The geometric multigrid (GMG) method can provide a solution in a way that is (nearly) optimal in terms of computational and memory complexity.  However GMG methods require a predetermined hierarchical structure,they are difficult to use in practice.  On the other hand, the algebraic multigrid (AMG) methods build the hierarchical structure by using only the information of the matrix and meanwhile maintain competitive performance for certain model problems.  In fact, for such problems, AMG methods are very attractive because they are highly efficient and easily to implement.  There are many different types of AMG methods: the classical AMG \cite{Stuben1983,Brandt1984}, smoothed aggregation AMG \cite{Vanek1996,Brezina2005} and AMGe \cite{Jones2002,Lashuk2008}, etc.  Given the need to solve extremely large systems, researchers have developed and implemented a number of parallel AMG methods for CPU. BoomerAMG (included in the Hypre package \cite{Falgout2002}) is the parallelization of the classical AMG methods and their variants, whereas others including ML \cite{Gee2006} focus on the parallel versions of smoothed aggregation AMG methods. In addition, some parallel implementations of the AMG methods are used in commercial software such as NaSt3DGPF \cite{nast} and SAMG \cite{samg}.

Not only are the researchers rapidly developing algorithms, they are doing the same with the hardware.  GPUs based on single instruction multiple thread (SIMT) hardware architecture have been provided an efficient platform for large-scale scientific computing since November 2006, when NVIDIA released the Compute Unified Device Architecture (CUDA) toolkit.  The CUDA toolkit made programming on GPU considerably easier than it had been previously in large part. 
Consequently, many sparse linear algebra and solver packages based on CUDA have been developed: MAGMA \cite{Tomov2010},CULA \cite{Humphrey2010}, the CUDPP library \cite{Sengupta2007}, NIVIDA CUSPARSE library \cite{Bell2009}, IBM SpMV library \cite{Baskaran2009}, Iterative CUDA \cite{Klockner}, Concurrent Number Cruncher \cite{Buatois2007} and CUSP \cite{Garland2010}). 

MG methods have also been parallelized and implemented on GPU in a number of studies.  GMG methods as the typical cases of AMG have been implemented on GPU firstly\cite{Goodnight2003,Bolz2003,Goddeke2008,Grossauer2008,Feng2010,chensong2012}.
These studies demonstrate that the speedup afforded by using GPUs can result in GMG methods achieving  a high level of performance on CUDA-enabled GPUs.  
However, to the best of our knowledge of the task, paralleling an AMG method on GPU or CPU remains very challenging mainly due to the sequential nature of the coarsening processes (setup phase) used in AMG methods.  In most  AMG algorithms, coarse-grid points or basis are selected sequentially using graph theoretical tools (such as maximal independent sets and graph partitioning algorithms) and the coarse-grid matrices are constructed by a triple-matrix multiplication.  Although extensive research has been devoted to improving the performance of parallel coarsening algorithms, leading to marked improvements on CPU  \cite{Cleary1998,Tuminaro2000,Krechel2001,Henson2002,DeSterck2006,Chow2006,Joubert2006} and on GPU \cite{Tuminaro2000,Bell2011,Kraus} over time, the setup phase is still considered the major bottleneck in parallel AMG methods.  On the other hand, the task of  designing an efficient and robust parallel smoother in the solver phase is no less  challenging. 


In this paper, we design a new parallel AMG method based on the unsmoothed aggregation AMG (UA-AMG) method together with what is known as the nonlinear AMLI-cycle (K-cycle) \cite{Kraus.J2002,Notay2008,Hu.X;Vasilevski.P;Xu.J2011a}. The UA-AMG method is attractive because of the simplicity of its setup procedure and its low computational cost.  We designed our new AMG method to overcome the difficulties of the setup phase, i.e. algebraic coarsening process. And we did this by using the geometric grid on the finest level in order to build the coarse levels instead of using the pure algebraic coarsening characteristic of most AMG methods. The main idea is to construct an auxiliary structured grid based on   information from the finest geometric grid, and select a simple and fixed coarsening algorithm that allows explicit control of the overall grid and operator complexities of the multilevel solver.  When an auxiliary structured grid is used,  the coarsening process of the UA-AMG method, i.e. the construction of aggregates, can easily be done by choosing the elements that intersect certain local patches of the auxiliary grid.  This auxiliary grid and its relationship with the original unstructured grid can be effectively managed by a quadtree in 2D (octree in 3D).  The special structure of the auxiliary grid narrows the bandwidth of the coarse grid matrix (9-point stencil in 2D and 27-point stencil in 3D). And, this narrowing gives us explicit control of the sparsity pattern of the coarse-grid matrices and reduces the operator complexity of the AMG methods. These features control the work per thread which is considered to be advantage for parallel computing, especially on GPU.  Moreover, due to the regular sparsity pattern,  we can use the ELLPACK format to significantly speed up the sparse matrix-vector multiplication on GPU  (see \cite{Bell2008} for a discussion of different sparse matrix storage formats).  As our new parallel AMG method is based on the UA-AMG framework, there is no need to form the prolongation and restriction matrices explicitly. Both of them can be done efficiently in parallel by using the quadtree of the auxiliary structured grid.  In addition, the auxiliary grid allows to use the colored Gauss-Seidel smoother without any coloring the matrices.  This improves not only the convergence rate but also the parallel performance of the UA-AMG method.

The remainder of the paper is organized as follows: In Section \ref{sec:UAAMG}, we review the  UA-AMG method.  Then in Section \ref{sec: SG-aggregation}, we introduce the algorithm we used to form the aggregations based on the geometric information.  Next, we demonstrate the parallelization of  our algorithm on GPU in Section \ref{sec:ParaGPU}. In Section \ref{sec:result}, we discuss the results of our numerical experiments in regard to  our proposed algorithm's performance on GPU. Finally, we summarize our findings with some concluding remarks in Section \ref{sec:conclusion}.

\section{Unsmoothed Aggregation AMG}\label{sec:UAAMG}

In this section, because our parallel AMG method is based on the UA-AMG method, we recall the unsmoothed aggregation AMG method \cite{CNM:CNM1640090804,Braess1995}.   The coarsening in the UA-AMG method is performed by simply aggregating the unknowns, and the prolongation and restriction are Boolean matrices that characterize the aggregates.  It has been shown that under certain conditions the Galerkin-type coarse-level matrix has a controllable sparsity pattern \cite{Kim2003a}.  These properties make the UA-AMG method more suitable for parallel computation than are traditional AMG methods, such as the classical AMG and  the smoothed aggregation AMG, especially on GPUs.  Recent work also shows that this approach is both theoretically well founded and practically efficient, providing that we form the aggregates in an appropriate way and apply certain enhanced cycles, such as the AMLI-cycle or the Nonlinear AMLI-cycle (K-cycle) \cite{Kraus.J2002,Notay2008,Hu.X;Vasilevski.P;Xu.J2011a}. 

Given the $k$-th-level matrix $A_k \in \mathbb{R}^{n_k \times n_k}$, in the UA-AMG method we define the prolongation matrix $P_{k-1}^{k}$ from a non-overlapping partition of the $n_k$ unknowns at level $k$ into the $n_{k-1}$ nonempty disjoint sets $G_j^{k}$, $j=1,\dots, n_{k-1}$, which are referred to as  aggregates. There are many different approaches to generating the aggregates. However, in most of these approaches, the aggregates are chosen one by one, i.e. sequentially.  And, standard parallel aggregation algorithms are variants of the parallel maximal independent set algorithm, which is choosing the points of the independent set based on the giving random numbers for every points.  Therefore, the quality of the aggregates cannot be guaranteed. 

In this paper, we build the aggregates based on the structured auxiliary gird which constructed by the geometric information of the original unstructured grid. All aggregates can be chosen independently. This makes the aggregation step suitable for parallel computing. The details of our new algorithm are presented in the next section.  Here, let us assume that the aggregates has been formed and are ready to use.  Once the aggregates are constructed, the prolongation $P_{k-1}^{k}$ is an $n_{k} \times n_{k-1}$ matrix given by 
\begin{equation} \label{def:P}
(P_{k-1}^{k})_{ij} = 
\begin{cases}
1 & \text{if} \ i \in G_j^{k}  \\
0 & \text{otherwise}
\end{cases} \quad
i = 1, \dots, n_k, \quad j = 1, \dots, n_{k-1}. 
\end{equation}
With such a piecewise constant prolongation $P_{k-1}^{k}$, the Galerkin-type coarse-level matrix $A_{k-1} \in \mathbb{R}^{n_{k-1} \times n_{k-1}}$ is defined by
\begin{equation} \label{def:Ac}
A_{k-1} = (P_{k-1}^{k})^t A_k (P_{k-1}^{k}).
\end{equation}
Note that the entries in the coarse-grid matrix $A_{k-1}$ can be obtained from a simple summation process:
\begin{equation}\label{eqn:acsum}
(A_{k-1})_{ij} = \sum_{k \in G_{i}} \sum_{l \in G_{j}} a_{kl}, \quad i,j = 1, 2, \cdots, n_{k-1}.
\end{equation}

The enhanced cycles are often used for the UA-AMG method in order to achieve uniform convergence for model problems such as the Poisson problem \cite{Notay2008,Notay2010,Hu.X;Vasilevski.P;Xu.J2011a} and M-matrix \cite{Kraus.J2002}.  In this paper, we use the nonlinear AMLI-cycle, also known as the variable AMLI-cycle or the K-cycle because it is parameter-free. The nonlinear AMLI-cycle uses the Krylov subspace iterative method to accelerate the coarse-grid correction. Therefore, for the completeness of our presentations, let us first recall the nonlinear preconditioned conjugate gradient (PCG) method originated in \cite{Fletcher1964}.  Algorithm~\ref{alg:nonlinear-CG} is a simplified version designed to address symmetric positive definite (SPD) problems.  The original version was meant for more general cases including nonsymmetric and possibly indefinite matrices.  Let~$\hat{B}_k[\cdot]: \mathbb{R}^{n_k} \rightarrow \mathbb{R}^{n_k}$ be a given nonlinear preconditioner intended to approximate the inverse of~$A_k$. We now formulate the nonlinear PCG method that can be used to provide an iterated approximate inverse to~$A_k$ based on the given nonlinear operator~$\hat{B}_k[\cdot]$. This procedure gives another nonlinear operator~$\tilde{B}_k[\cdot]: \mathbb{R}^{n_k} \rightarrow \mathbb{R}^{n_k}$, which can be viewed as an improved approximation of the inverse of~$A_k$.

\begin{algorithm}[htbp!]
\caption{Nonlinear PCG Method} \label{alg:nonlinear-CG}
\begin{flushleft}
Assume we are given a nonlinear operator~$\hat{B}_k[\cdot]$~to use as a preconditioner.
Then, for~$\forall f \in \mathbb{R}^{n_k}$, $\tilde{B}_k[f]$~is defined as follows:
\begin{enumerate}
\item [] {\bf Step 1.} Let~$u_0 = 0$~and~$r_0 = f$. Compute~$p_0 = \hat{B}_k[r_0]$. Then let
\[
u_1 = \alpha_0 p_0, \ \text{and} \ r_1=r_0 - \alpha_0A_kp_0, \ \text{where} \ \alpha_0 = \frac{(r_0, p_0)}{(p_0,p_0)_{A_k}}.
\] 
\item [] {\bf Step 2.} For~$i=1,2,\cdots,n-1$, compute the next conjugate direction
\begin{equation*}\label{eqn:conjugate_direction}
p_i = \hat{B}_k[r_i] + \sum_{j=0}^{i-1} \beta_{i,j} p_j, \ \text{where} \ \beta_{i,j} = -\frac{(\hat{B}_k[r_i], p_j)_{A_k}}{(p_j,p_j)_{A_k}}.
\end{equation*}
Then the next iteration is
$u_{i+1} = u_{i} + \alpha_i p_i, \ \text{where} \ \alpha_i = \frac{(r_i, p_i)}{(p_i,p_i)_{A_k}}$, 
and the corresponding residual is
$r_{i+1} = r_{i} - \alpha_i A_k p_i$.

\item [] {\bf Step 3.} Let~$\tilde{B}_k[f] := u_n$.
\end{enumerate}
\end{flushleft}
\end{algorithm}

We also need to introduce a smoother operator~$R_k: \mathbb{R}^{n_k} \rightarrow \mathbb{R}^{n_k}$~in order to define the multigrid method.  In general, all smoothers, such as the Jacobi smoother and the Gauss-Seidel (GS) smoother, can be used here.  In this paper, we use parallel smoothers which is based on auxiliary grid. In subsection \ref{sec:parallel_smoother}, we discuss the algorithm and implementation of the smoother in details.  Here, it is assumed that we already have a smoother.  Now, thanks to Algorithm~\ref{alg:nonlinear-CG}, we can recursively define the nonlinear AMLI-cycle MG operator as an approximation of~$A_k^{-1}$ (see Algorithm~\ref{alg:Kcycle}).

\begin{algorithm}
\caption{Nonlinear AMLI-cycle MG:~$\hat{B}_{k}[\cdot]$} \label{alg:Kcycle}
\begin{flushleft}
Assume that $\hat{B}_1[f]=A_1^{-1}f$, and $\hat{B}_{k-1}[\cdot]$ has been defined, then for $f \in \mathbb{R}^{n_k}$
\end{flushleft}
\begin{enumerate}
\item [] {\bf Pre-smoothing} $u_1 =  R_{k}f$
\item [] {\bf Coarse-grid correction} $u_2 = u_1 +  (P_{k-1}^{k})^{t}  \tilde{B}_{k-1}[P_{k-1}^{k}(f-A_k u_1)]$, where~$\tilde{B}_{k-1}$~is implemented as in Algorithm~\ref{alg:nonlinear-CG} with~$\hat{B}_{k-1}$~as the preconditioner;
\item [] {\bf Post-smoothing} $\hat{B}_k[f] := u_2 + R_k^t(f - A_k u_2)$.
\end{enumerate}
\end{algorithm}

Our parallel AMG method is mainly based on Algorithm~\ref{alg:Kcycle} with prolongation $P_{k-1}^{k}$ and coarse-grid matrix $A_{k-1}$ defined by \eqref{def:P} and \eqref{def:Ac} respectively.  The main idea of our new parallel AMG method is to use an auxiliary structured grid to (1) efficiently construct the aggregate in parallel; (2) simplify the construction of the prolongation and the coarse-grid matrix; (3) develop a robust and effective colored smoother; and (4) better control the sparsity pattern of the coarse-grid matrix and working load balance.  

\section{Auxiliary Grid Aggregation Method} \label{sec: SG-aggregation}
In this section, we discuss how to form the aggregates on each level, which is the essential procedure in the setup phase of the UA-AMG method.  Aggregates are usually generated by finding the maximal independent set of the corresponding graph of the matrix.  This procedure is mainly sequential and difficult to implement in parallel.  Moreover, the shape of the coarse grid and of the sparsity pattern of the coarse-grid matrix cannot usually be controlled, which increases the computational complexity of the AMG methods.  This lack of control and corresponding complexity makes the traditional aggregation approaches less favorable for parallel computation than other AMG methods, especially for GPUs. However, our main focus here is discretized partial differential equations, for which detailed information on the underlying geometric grid is generally available, although we usually only have  access to the finest-level unstructured grid.  We propose an aggregation method wherein information from the finest grid is used to select a simple and fixed coarsening that allows the overall grid and the operator complexities to be explicitly controlled.  We use the geometric information of the underlying unstructured grid to construct an auxiliary structured grid and build a hierarchical structure on the auxiliary structured grid.  The aggregates on each level are formed based on the hierarchical structure of the structured grid.  Moreover, the auxiliary structured grid and its hierarchical structure can be managed effectively by a quadtree in 2D (or an octree in 3D).  Another main advantage of the auxiliary grid is that we can easily color the grids on every level, which makes using the colored Gauss-Seidel smoother feasible. A colored Gauss-Seidel smoother helps both to improve the overall convergence behavior of our AMG method and to achieve a high level of parallelism.   Details regarding the colored GS smoother are addressed in subsection~\ref{sec:parallel_smoother} .  In this section, we describe the auxiliary grid aggregation method for the 2D case.  However, we would like to point out that it can easily be generalized to the 3D case. 

\subsection{Auxiliary Structured Grid and Quadtree}
Assume that the computational domain is $\Omega$, which is triangulated into shape-regular elements $\tau_i$ such that 
\[ 
\bigcup_{i} \closure{\tau}_i = \closure[2]{\Omega}\in \mathbb{R}^d, \quad d = 2 \ \text{or} \ 3.  
\]
We denote the unstructured grid by $\mathcal{T} = \{  \tau_i \}$.  


In order to construct an auxiliary structured grid and build its hierarchical structure, we start with another domain $\Omega^0_0$ such that $\closure{\Omega}\subset\closure{\Omega}^0_0$.  A simple way to choose the domain $\Omega^0_0$ is to find a rectangular domain that covers the whole domain $\Omega$. i.e,
\begin{equation}\label{eqn:auxdomain}
\Omega^0_0 = (a_1,b_1)\times(a_2,b_2)
\end{equation}
where
\begin{equation} \label{def:ab12}
 a_1=\min{x},\quad a_2 = \min{y},\quad b_1 = \max{x},\quad b_2=\max{y}, \quad \text{for} \ (x,y) \in \closure{\Omega}.
\end{equation}

The way to construct the hierarchical structure is by constructing a region quadtree~\cite{Finkel1974}. The region quadtree represents a partition of space in two dimensions by breaking the region into four quadrants, subquadrants, and so on, with each leaf node containing data corresponding to a specific subregion \cite{DeBerg2000}. The region tree $T$ has the following properties:
\begin{enumerate}
\item A subregion $\Omega^{k}_{j}$ is a leaf when either $k = L$ or the subregion has 4 children and $k < L$ (where L is the depth of the tree);  
\item A subregion $\Omega^{k}_{j}$, $k < L$ is defined as the union of its $4$ children:
\begin{equation}\label{eqn:quadtree}
\Omega^{k}_{j} = \bigcup_{i} \Omega^{k+1}_{i}.
\end{equation}
\item  For all subregions on the same level $k$, 
\[\bigcup_{j} \Omega^{k}_{j} = \Omega.\]
\end{enumerate}

\begin{figure}[!htp]
\small
\centering
\includegraphics[width=0.45\textwidth]{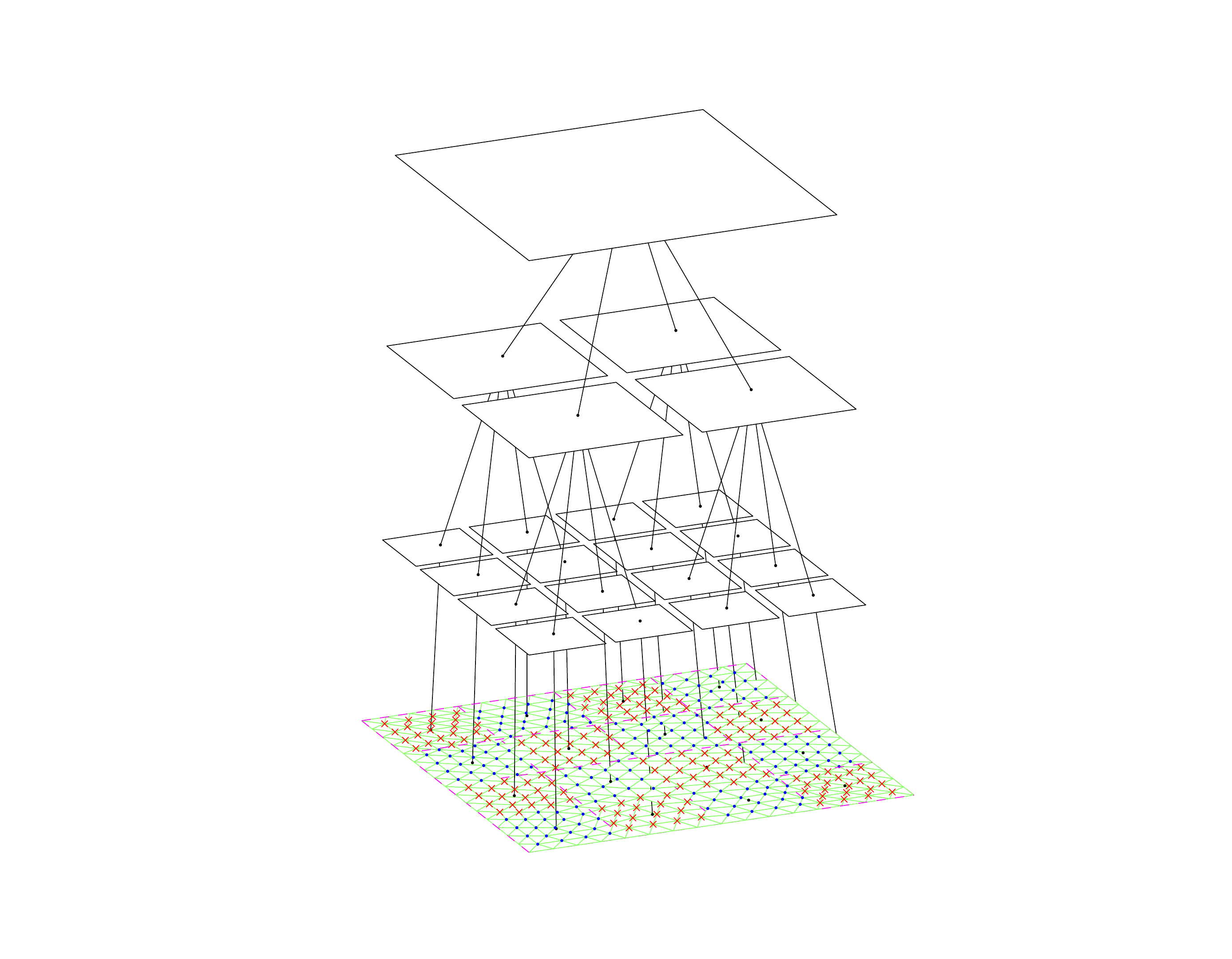}
\includegraphics[width=0.45\textwidth]{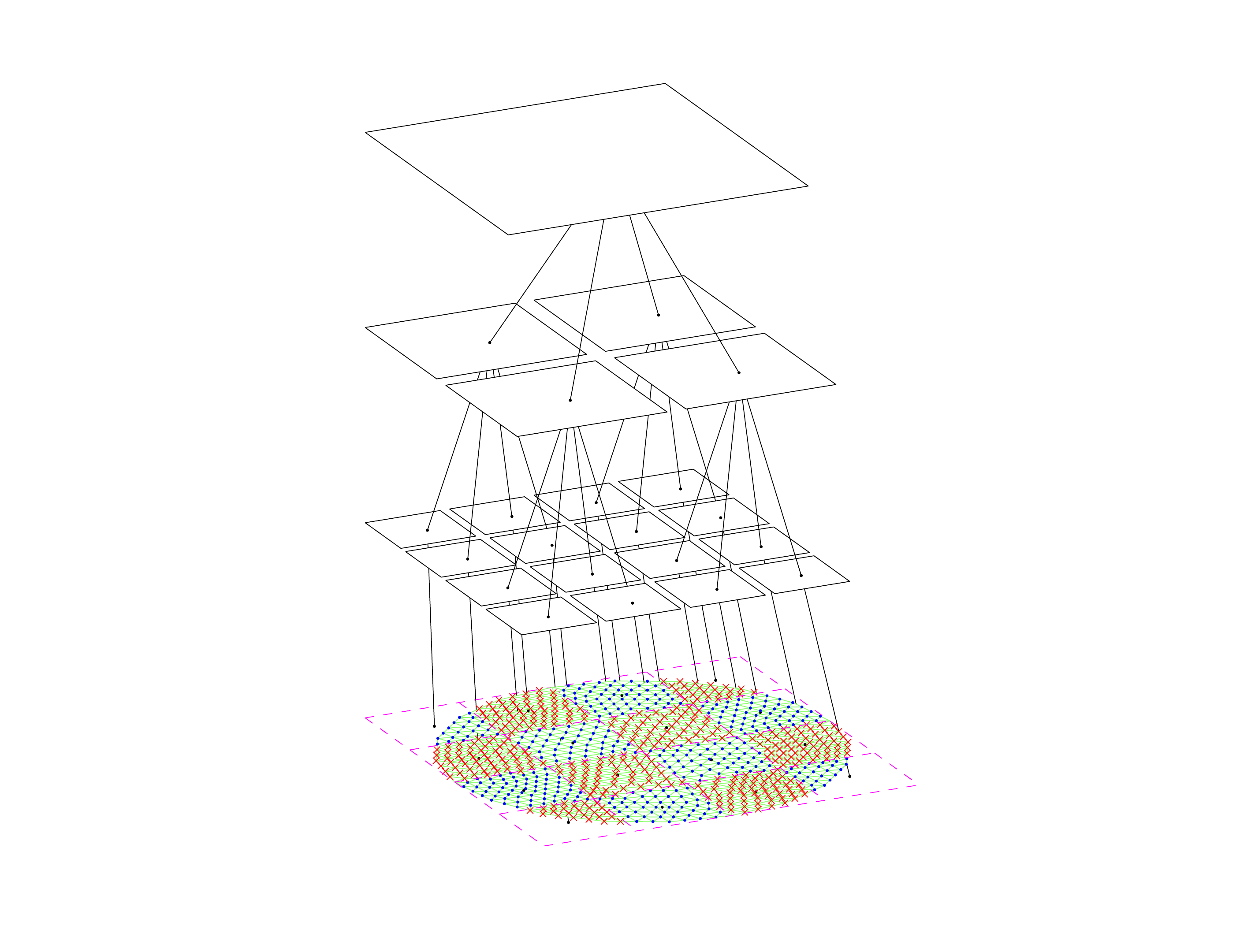}
\caption{Examples of the region quadtree on different domains.}
\label{fig:quadtree}
\end{figure}

Now we construct the region quadtree.  We use the domain $\Omega^0_0$ as a root and divide it into $4$ children subregions, such that 
\[
\begin{matrix}
 \Omega^1_0= (a_1,b_1')\times(a_2,b_2'),\quad& \Omega^1_1= (a_1',b_1)\times(a_2,b_2')\\
  \Omega^1_2= (a_1,b_1')\times(a_2',b_2),\quad& \Omega^1_3= (a_1',b_1')\times(a_2',b_2)
\end{matrix}
\]
where $a_1'=b_1'=(a_1+b_1)/2$ and $a_2'=b_2'=(a_2+b_2)/2$. Next, we apply the same process to the $4$ children $\Omega^1_0,\Omega^1_1,\Omega^1_2,\Omega^1_3$ and then to their children.  By doing this recursively, we can generate a region quadtree $T$ with root $\Omega^0_0$.  Because $T$ is a full quadtree, it is easy to check that, for each subregion $\Omega_i^k$ on level $k$ for $i = t_1 + 2^k t_2$, $0 \leq t_1, t_2 \leq 2^k-1$, we have
\begin{equation}\label{def:subregion}
\Omega^k_i = (a_1+t_1\frac{b_1-a_1}{2^k}, a_1+(t_1+1)\frac{b_1-a_1}{2^k})\times (a_2+t_2\frac{b_2-a_2}{2^k}, a_2+(t_2+1)\frac{b_2-a_2}{2^k}). 
\end{equation}
Figure \ref{fig:quadtree} gives two examples of the region quadtrees on a unit square domain and a circle domain.

It is easy to see that on level $k$, there are $4^k$ subregions, and they form a structured grid.  Now, we determine the depth $L$ of the region quadtree, and the subregion on that level $\{ \Omega^L_i \}$ is our auxiliary structured grid.  In order to limit the computational complexity and  ensure that the aggregates on the fine level do not become too small, we choose the level $L$, such that $4^L < N$, where $N$ is the number of vertices of the fine-level unstructured grid.  We can choose $L = \lfloor\frac{\log(N)}{\log4}\rfloor$, where the notation $n = \lfloor m \rfloor$ means that n is the biggest possible integer such that $n \leq m$.  




If $d = 3$, we can generate an octree as the auxiliary structured grid in the same way. Therefore, choose $L = \lfloor\frac{\log(N)}{\log(8)}\rfloor$.
\subsection{Aggregation Based on the Auxiliary Grid}
After we have the auxiliary structured grid and the region quadtree, we can form the aggregates of all the levels.  Here, we need to distinguish between two different cases (1) level $L$ and (2) all the coarse levels $0 < k < L$.
\begin{itemize}
\item {\bf Level $L$.} On this level, we need to aggregate the degree of freedoms (DoFs) on the unstructured grid. Basically, all the DoFs in the subregion $\Omega_i^L$ form one aggregate $i$.  Therefore, there are $4^L$ aggregates on this level ($8^L$ in 3D).  Let us denote the coordinate of DoF $j$ on the unstructured grid by $(x_j, y_j)$, such that the aggregate $G_i^L$ is defined by
\begin{equation} \label{def:agg_L} 
G_i^L = \{ j: (x_j, y_j) \in \Omega_i^L \}. 
\end{equation} 
Using~equation (\ref{def:subregion}), we can see that the aggregation generating procedure can be done by checking the coordinates of the DoFs. 
\item {\bf Coarse level $0<k<L$.} Start from the second level $L-1$ and proceed through the levels to the coarsest level $0$. Note that each level $k$ is a structured grid formed by all the subregions $\{ \Omega^k_i \}$, and note, too, that a DoF on the coarse level corresponds to a subregion on the coarse level. Therefore, we can use a fixed and simple coarsening whereby all the children of a subregion form an aggregate, i.e,
\begin{equation} \label{def:agg_k}
G_i^k = \{j: \Omega_j^{k+1} \ \text{is a child of} \ \Omega^k_i\}.
\end{equation}
Thanks to the quadtree structure, the children on the quadtree can be found by checking their indices. 
\end{itemize}

\begin{figure}[!htp]
\small
\centering
\includegraphics[width=0.48\textwidth]{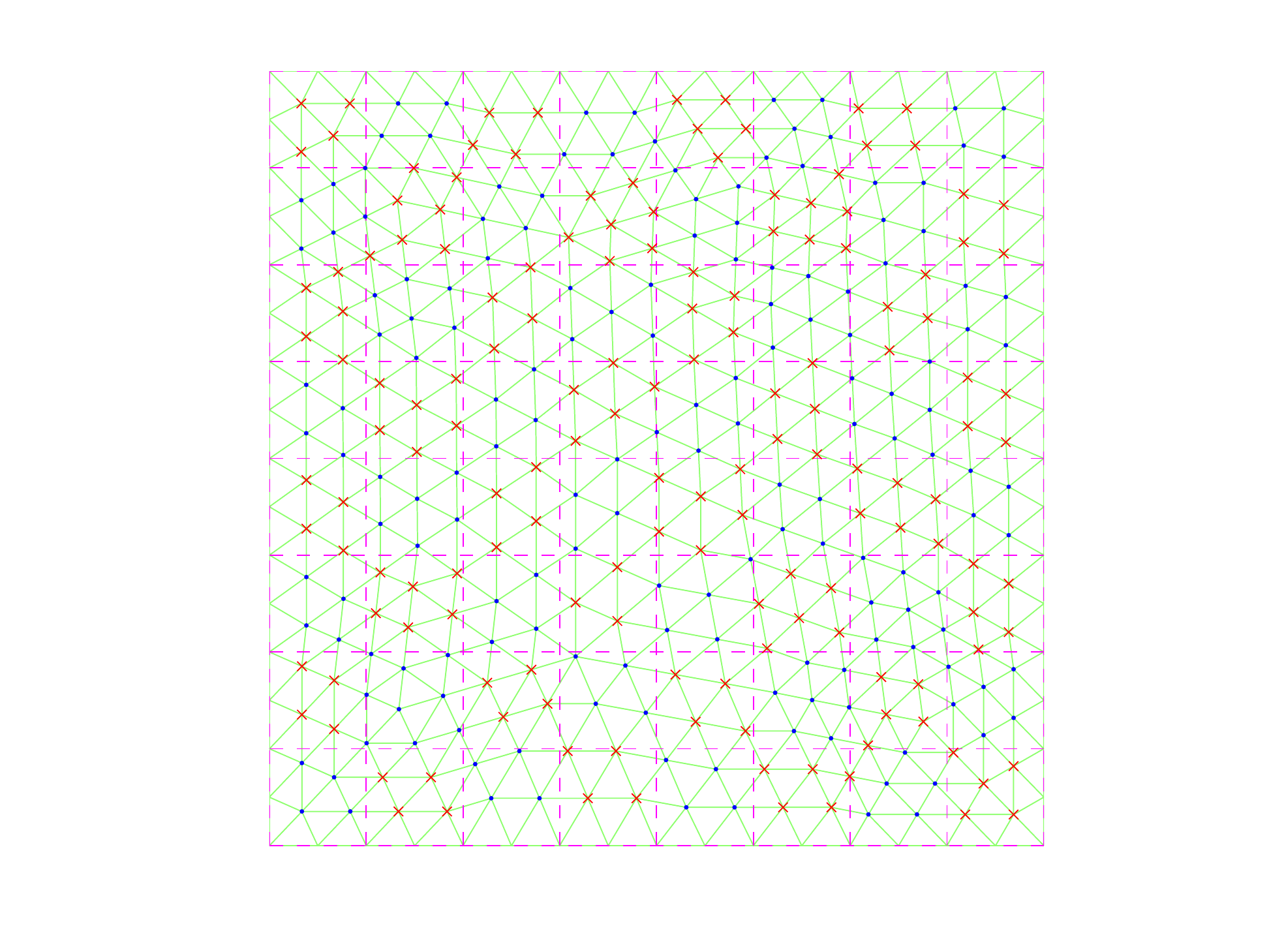}
\includegraphics[width=0.48\textwidth]{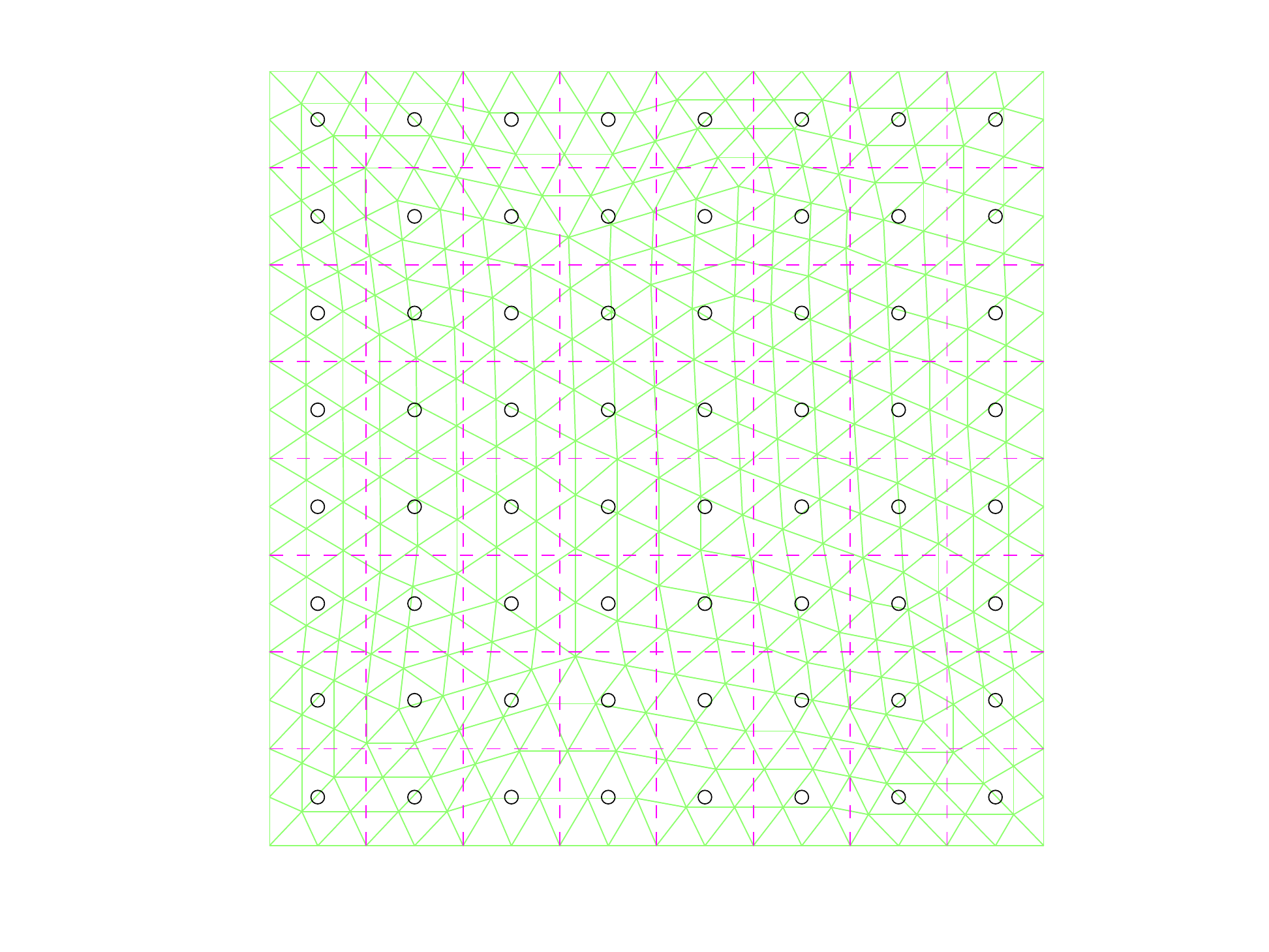}
\caption{Aggregation on level $L$.}
\label{fig:aggregation1}
\end{figure}

\begin{figure}[!htp]
\small
\centering
\includegraphics[width=0.49\textwidth]{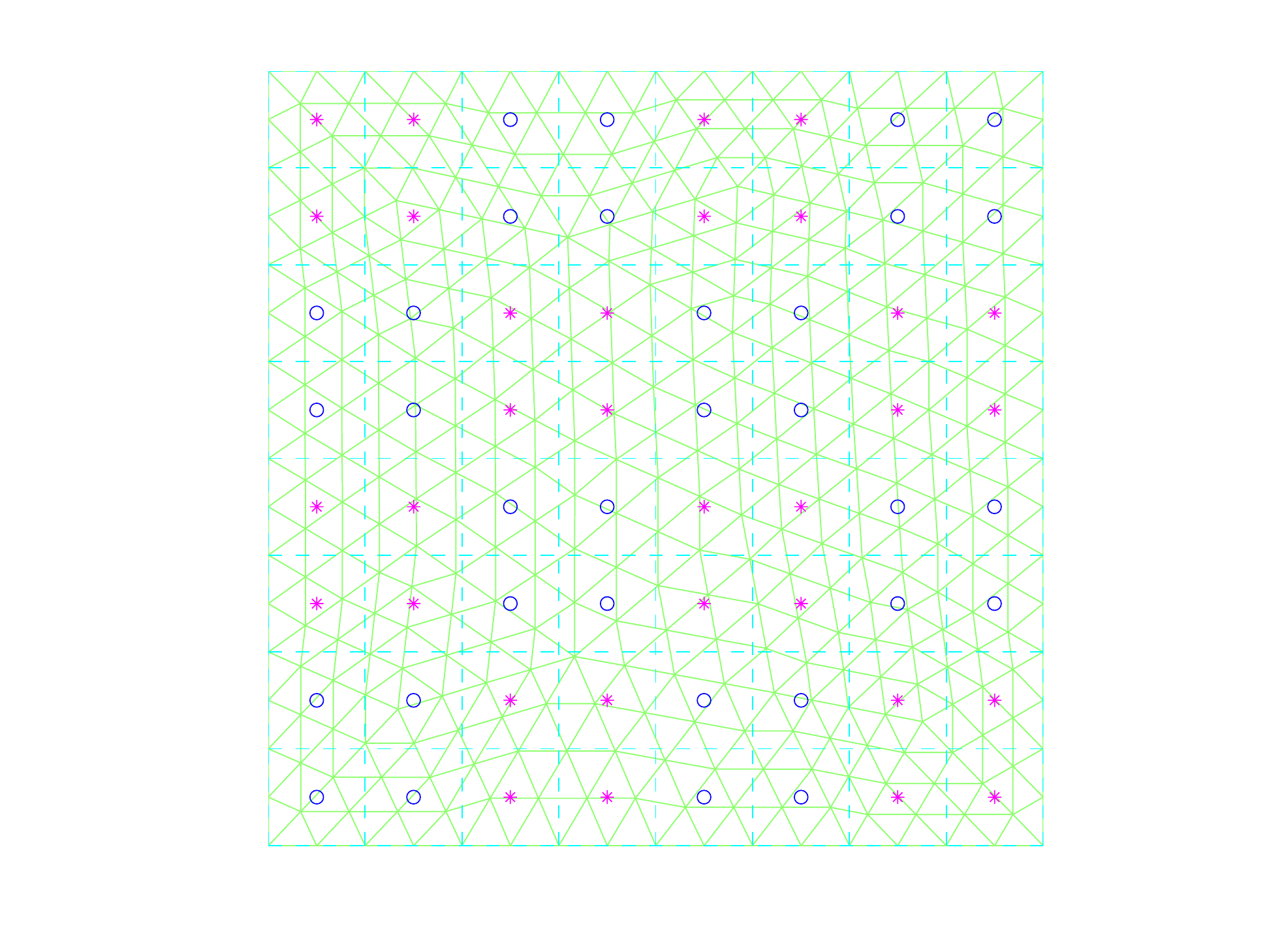}
\includegraphics[width=0.41\textwidth]{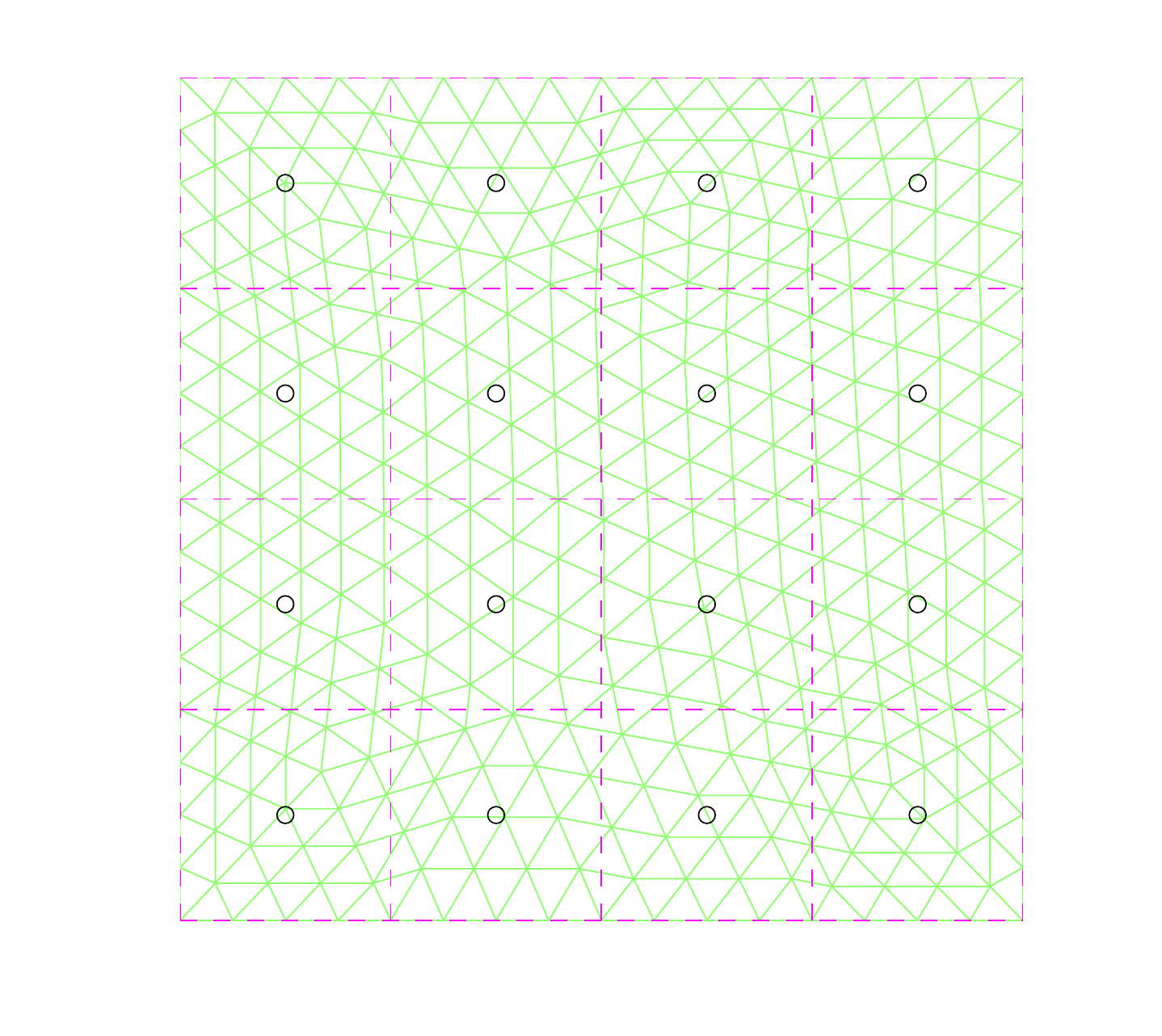}
\caption{Aggregation on the coarse levels.}
\label{fig:aggregation2}
\end{figure}

Figures~\ref{fig:aggregation1}~and~\ref{fig:aggregation2}~ show the aggregation on level $L$ and on the coarse levels respectively.  The detailed parallel algorithms and implementation of them are discussed in the next section.



\section{Parallel Auxiliary Grid AMG Method on GPU}\label{sec:ParaGPU}
Many research groups have shown that GPU can accelerate multigrid methods, making them perform faster and better than on CPU (see \cite{Goodnight2003,Bolz2003,Goddeke2008,Grossauer2008,Feng2010,Goddeke2008,Haase2010}). Well-implemented GPU kernels can be more than 10 times faster than their CPU counterparts.  

In this section, we discuss the parallelization of our auxiliary grid AMG method on GPU, as follows:
\begin{enumerate}
\item  \emph{Sparse matrix-vector multiplication:} just as many other iterative methods do,  sparse matrix-vector multiplication (SpMv) plays an important rule in our AMG method. Moreover, because we apply the nonlinear AMLI-cycle and use our AMG method as a preconditioner for the Kyrlov iterative method,  we need an efficient parallel SpMV implementation on GPU. This means we should use a sparse matrix storage format that has a regular memory access pattern when we do the SpMV. 
\item \emph{Parallel aggregation algorithm:} in the coarsening procedure, aggregation is one of the most time-consuming steps.  In the previous section, we introduced our auxiliary-grid-based aggregation method and show its potential for efficient parallelization. We discuss the parallelization of the auxiliary-grid aggregation algorithm in this section.
\item \emph{Prolongation and restriction:}  an advantage of the UA-AMG method is that there is no need to explicitly form prolongation $P$ or restriction $R$, which usually is the transpose of $P$. This is because both P and R are just Boolean matrices that characterize the aggregates.  However, we do need their actions in the AMG method to transfer the solution and the residual between levels. And, we also need to be able to transfer the solution and the residual between levels in parallel.
\item \emph{Coarse-level matrix:} the coarse-level matrix is usually computed by a triple-matrix multiplication~\eqref{def:Ac}. This computation is considered to be the most difficult part to parallelize in the AMG method on GPU, as speedup can be as low as $1.2$ when compared with its CPU implementation (see \cite{Bell2011}).  However, thanks to the auxiliary structured grid, we have fixed sparsity pattern on each coarse level, which makes the triple-matrix multiplication much easier to compute. 
\item \emph{Parallel smoother:} the smoother is the central component of the AMG method. It is usually a linear iteration. For the parallel AMG method, however, it is difficult to design an efficient and parallel smoother.  Generally speaking, smoothers that have a good smoothing property, such as the Gauss-Seidel method, are difficult to parallelize, whereas, smoothers, like the Jacobi method, that are easy to parallelize cannot usually achieve a good smoothing property.  In our parallel AMG method, the auxiliary structured grid make it possible to use a colored Gauss-Seidel method, which maintains a good smoothing property and is easy to parallelize. 
\end{enumerate}

By combining these components with the nonlinear AMLI-cycle, we have established our parallel auxiliary grid AMG method, which is presented at the end of this section.


\subsection{Sparse Matrix-Vector Multiplication on GPUs}
An efficient SpMV algorithm on GPU requires a suitable sparse matrix storage format.  How different storage formats perform in SpMV is extensively studied in \cite{Bell2008}. 
This study shows that the need for coalesce accessing of the memory makes matrix formats such as compressed row storage (CSR) format and compressed column storage (CSC) format widely used for the iterative linear solvers on CPU but inefficient on GPU. According to \cite{Bell2008}, when each row has roughly the same nonzeros, ELLPACK (ELL) storage format is one of the most efficient sparse matrix storage formats on GPU.  

Let us recall the ELL format, an $M \times N$ sparse matrix with at most $K$ nonzeros per row is stored in two $M \times K$ arrays: (1) array \texttt{Ax} stores nonzero entries, and (2) array \texttt{Aj} stores the column indices of the corresponding entries.  Both arrays are stored in column-major order to both improve the efficiency when accessing the memory on GPU~\cite{Bell2009}, see Figure \ref{fig:ellform}.

\begin{figure}[!htp]
\small
\centering
\begin{equation*}
A = \begin{pmatrix}
4&-2&0&0\\
-1&2&-1&0\\
0& 0& -3&4
\end{pmatrix},
\quad  \texttt{Aj}=
\begin{bmatrix}
0&1&-1\\
0&1&2\\
2& 3&-1
\end{bmatrix}
\quad 
\texttt{Ax}=
\begin{bmatrix}
4&-2&0\\
-1&2&-1\\
 -3&4&0
\end{bmatrix}
\end{equation*}
\flushleft{\textbf{Column-major ordering of \texttt{Aj} and \texttt{Ax}: }}
\begin{tabular}{rccccccccc}
\text{Aj: }&[0&0& 2& 1& 1& 3& -1& 2& -1]\\
\text{Ax: }&[4& -1& -3& -2& 2& 4& 0& -1& 0]\\
\end{tabular}
\flushleft{\textbf{Memory access of SpMV: one thread handle one row}}
\begin{align*}
\texttt{Ax} &= \left[\begin{array}{ccccccccc} 4 & -1 & -3 & -2 & 2 & 4 & 0& -1& 0  \end{array}\right] \\  
\text{Step 1} &= \left[\begin{array}{ccccccccc}  0 &  1 & 2 & \ \ & \ \ & \ \ & \ \ & \ \  & \ \
  \end{array}\right] \\
\text{Step 2} &= \left[\begin{array}{ccccccccc} \ \ & \ \  & \ \   & 0 & 1 & 2 & \ \ & \ \ & \ \ 
 \end{array}\right]  \\
\text{Step 3} &= \left[\begin{array}{ccccccccc} \ \  & \ \  & \ \ & \ \ & \ \ & \ \ & \ \  & 2 & \   
\end{array}\right]
\end{align*}
\caption{Sparse matrix representation using the ELL format and the memory access pattern of SpMv.}
\label{fig:ellform}
\end{figure}

It is not always possible to use the ELL format for AMG method, because the sparsity patterns of the coarse-level matrices can become denser than the fine-level matrices and the number of nonzeros per row can vary dramatically.  However, the ELL format is the perfect choice for our AMG method.  This is because the auxiliary structured grid and the region quadtree.  All the coarse-level matrices generated from the auxiliary grid have a $9$-point stencil structure and their sparsity patterns are fixed.  Based on these facts, we choose the ELL format to maximize the efficiency of SpMV. Moreover, we store the diagonal entries in the first column of \texttt{Aj} and \texttt{Ax} to avoid having to execute a searching process during smooth process and, therefore, further improve the overall efficiency of our parallel AMG method.



\subsection{Parallel Auxiliary Grid Aggregation}
As shown in Section \ref{sec: SG-aggregation}, we use an auxiliary structured grid and its hierarchical structure to form the aggregates.  A region quadtree is used to handle the hierarchical structure.  Here, we discuss how to efficiently parallelize this aggregation method.  Again, we distinguish two between cases: level $L$ and coarse levels $0 < k <L$.
\begin{itemize}
\item {\bf Level L:}  The aggregate $\{ G_i^L \}$ is defined by~\eqref{def:agg_L}, which means that we need to check the coordinates of each DoF, and determine the aggregates to which each  belongs.  Obviously, checking the coordinates of each DoF is a completely separate process; therefore, it can be efficiently performed parallel with other processes.  Algorithmically, we assign each thread to DoF $j$ and determine a subregion $\Omega_i^L$ such that the coordinates $(x_j, y_j) \in \Omega_i^L$, then $j \in G_i^L$.  Assume that $\Omega_i^L$ is labeled lexicographically, and then given $(x_j, y_j)$, the index $i$ is determined by 
\begin{equation*}
i =\lfloor \frac{y_j - a_2}{b_2 - a_2} \times 2^L \rfloor \times 2^L + \lfloor \frac{x_j - a_1}{b_1 - a_1} \times 2^L \rfloor,
\end{equation*}
where $a_1$, $a_2$, $b_1$, and $b_2$ are given by~\eqref{def:ab12}.  The output this parallel subroutine is an array \texttt{aggregation} that contains the information about the aggregates.  We would like to point out that because it substitutes for prolongation and restriction, \texttt{aggregation} plays an important rule in our AMG method.  Following is the psuedocode implemented using CUDA.

\begin{lstlisting}[title={Parallel aggregation on level $L$},label=agg_L]
__global__ void aggregation(int* aggregation, double* x, double* y, double xmax, double ymax, double xmin, double ymin, int L)
{
	/*  get thread ID */
	const unsigned int idx = blockIdx.x*blockDim.x+threadIdx.x;   
	const unsigned int powL = (int)pow(2.0, L);   
	/* check x coordinate  */
	const unsigned int  xIdx = (int)((x - xmin)/(xmax-xmin)*powL);   
	/* check y  coordinate */
	const unsigned int  yIdx = (int)((y - ymin)/(ymax-ymin)*powL);   
	/* label the aggregate */
	aggregation[idx] = yIdx*powL + xIdx;   
}
\end{lstlisting}

\item {\bf Coarse level $0 < k < L$:} The aggregate $\{G_i^k\}$ is defined by~\eqref{def:agg_k}. Therefore, we need to determine the root of each subregion $\Omega_j^{k+1}$ --  a task that can be performed in a completely parallel way.  Each thread is assigned to a DoF on the coarse level $k+1$, or equivalently to say a subregion $\Omega_j^k+1$, then compute the index $i$ on each thread in parallel by
\begin{equation*}
i = \lfloor \frac{j\%2^k}{2} \rfloor \times \frac{2^k}{2} +  \lfloor \frac{j/2^k}{2} \rfloor.
\end{equation*}
Following is the psuedocode implemented using CUDA.

\begin{lstlisting}[title={Parallel aggregation on coarse level $k$},label=agg_k]
__global__ void aggregation(int* aggregation, int k)
{
	/* get thread ID */
	const unsigned int idx = blockIdx.x*blockDim.x+threadIdx.x;   
	const unsigned int powk = (int)pow(2.0, k);   
	/* check x index */
	const unsigned int  xIdx = (int)((idx%powk)/2);
	/* check y index  */    
	const unsigned int  yIdx = (int)((idx/powk)/2);
	/*  label the aggregate */    
	aggregation[idx] = yIdx*(powk/2) + xIdx;   
}
\end{lstlisting}

\end{itemize}

\subsection{Parallel Prolongation and Restriction}
As noted in Section \ref{sec:UAAMG}, the prolongation and restriction matrices are piecewise constant and characterize the aggregates. Therefore, we do not form them explicitly in UA-AMG method.
\begin{itemize}
\item {\bf Prolongation:} Let $v^{k-1}\in \mathbb{R}^{n_{k-1}}$, so that the action $ v^k = P_{k-1}^{k}v^{k-1}$ can be written component-wise as follows:
\[
(v^k)_i = (P_{k-1}^{k}v^{k-1} )_{i} = (v^{k-1})_{j}, \quad i \in G^{k-1}_{j}
\]
Assign each thread to one element of $v^k$, and use the array \texttt{aggregation} to obtain information about $j \in G^{k-1}_i$, i.e., $i = \texttt{aggregation}[j]$, so that prolongation can be efficiently implemented in parallel. 
\item {\bf Restriction:} Let $v^{k}\in \mathbb{R}^{n_{k}}$, so that the action $(P_{k-1}^{k})^Tv^{k}$ can be written component-wise as follows:
\[
(v^{k-1})_i = ((P_{k-1}^{k})^Tv^{k} )_{i} = \sum_{j \in G^{k-1}_{i}}(v^{k})_{j}. 
\]
Moreover, when $0 < k < L$, assume that $i=t_1+2^{k-1}t_2$. The above formula can be expressed explicitly as
\[(v^{k-1} )_{i} = (v^{k})_{2t_1+2^k(2t_2)}+(v^{k})_{2t_1+1 +2^k(2t_2)}+(v^{k})_{(2t_1+1) +2^k(2t_2+1)}+(v^{k})_{2t_1 +2^k(2t_2+1)}.
\]
Therefore, each thread is assigned to an element of $v^{k-1}$, and use the array \texttt{aggregation} to obtain information about $j \in G^{k-1}_i$, i.e., find all $j$ such that $\texttt{aggregation}[j] = i$, on level $L$, or use the explicit expression on the coarse levels, such that the action of restriction can also be implemented in parallel.
\end{itemize}

\subsection{Parallel Construction of Coarse-level Matrices}
Usually, the coarse-level matrix is constructed by the triple-matrix multiplication \eqref{def:Ac}.  However, the triple-matrix multiplication is a major bottleneck of parallel AMG method on GPU \cite{Bell2011},  not only in terms of parallel efficiency, but also in regard to memory usage.  In general, the multiplication of sparse matrices comprises two steps.  First step is to determine the sparsity pattern by doing the multiplication symbolically, and the second step is to actually compute the nonzero entries.  In our AMG method, because of the presence of the auxiliary structured grid, the sparsity pattern is known and fixed. There is not necessary to perform the symbolic multiplication called for in the first step.  For the second step, according to \eqref{eqn:acsum}, the construction of coarse-level matrices in the UA-AMG method only involves simple summations over aggregates and can be performed simultaneously for each nonzero entry.  

We use the ELL format to store the sparse matrices. Therefore, we need to generate the index array $\texttt{Aj}_k$ and the nonzero array $\texttt{Ax}_k$ for a coarse-level matrix $A_{k}$.  
\begin{itemize}
\item {\bf Form $\texttt{Aj}_k$:}  due to the special structure of the region quadtree, all the coarse-level matrices have a 9-point stencil structure.  Based on this fact, $\texttt{Aj}_k$ is predetermined and can be generated in parallel without any symbolic matrix multiplication.  For a DoF, $i=t_1+2^kt_2$, on level k $(k \leq L)$, its $8$ neighbors are
\[
\begin{tabular}{lll}
$i_1=(t_1+1)+2^kt_2$,&$i_2=(t_1+1)+2^k(t_2+1)$,&$i_3=t_1+2^k(t_2+1)$,\\
$i_4=(t_1-1)+2^k(t_2+1)$,&$i_5=(t_1-1)+2^kt_2$,&$i_6=(t_1-1)+2^k(t_2-1)$,\\
$i_7=t_1+2^k(t_2-1)$,&$i_8=(t_1+1)+2^k(t_2-1)$.&
\end{tabular}
\]
And, the $i$-th row of the index array $\texttt{Aj}_k$ is, in Matlab notation, $\texttt{Aj}_k(i,:) = [i, i_1, i_2, \cdots, i_8]$, with some simple modifications for the DoFs on the boundaries.  It is easy to see that each thread can handle one entry in the array $\texttt{Aj}_k$ and complete it in in parallel. 
\item {\bf Form $\texttt{Ax}_k$:} The summation~\eqref{eqn:acsum} in the ELL format is as follows,
\begin{equation}\label{eqn:c-sum}
\texttt{Ax}_{k}(r,t) = \sum_{i \in G^{k}_{r}} \sum_{j \in G^{k}_q, \ q = \texttt{Aj}_{k}(r,t)} \texttt{Ax}_{k+1}(i,j), \quad r = 1, 2, \cdots, n_{k}, t = 0,1,\cdots,8.
\end{equation}
\end{itemize}
We use the Matlab notation again here.  Similar to the action of restriction, the searching for  $i \in G^{k}_{r}$ and $j \in G^{k}_q, \ q = \texttt{Aj}_{k}(r,t)$ can be done in parallel with the help of \texttt{aggregation}.  Together with the parallel summation, the array $\texttt{Ax}_k$ can be formed in parallel.  Moreover, on the coarse level $k<L$, we can write the summation explicitly as the action of restriction.  

\subsection{Parallel Smoothers Based on the Auxiliary Grid} \label{sec:parallel_smoother}
An efficient parallel smoother is crucial for the parallel AMG method.  For the sequential AMG method, Gauss-Seidel relaxation is widely used and has been shown to have a good smoothing property.  However the standard Gauss-Seidel is a sequential procedure that does not allow efficient parallel implementation. To improve the arithmetic intensity of the AMG method and to render it more efficient for parallel computing, we introduce colored Gauss-Seidel relaxation.  For example, on a structured grid with a $5$-point stencil (2D) or a $7$-point stencil (3D), the red-black Gauss-Seidel smoother is widely used for parallel computing because it can be efficiently parallelized (it is Jacobi method for each color) and still maintain a good smoothing property. In fact, in such a way, the red-black Gauss-Seidel smoother works even better than standard Gauss-Seidel smoother for certain model problems.  However, for an unstructured grid, coloring the grid is challenging, and the number of colors can be dramatically high for certain types of grid.  Therefore, how to apply a colored Gauss-Seidel smoother is often unclear, and  sometimes it may even be impossible to do so.  

\begin{figure}[!htp]
\small
\centering
\includegraphics[width=0.7\textwidth]{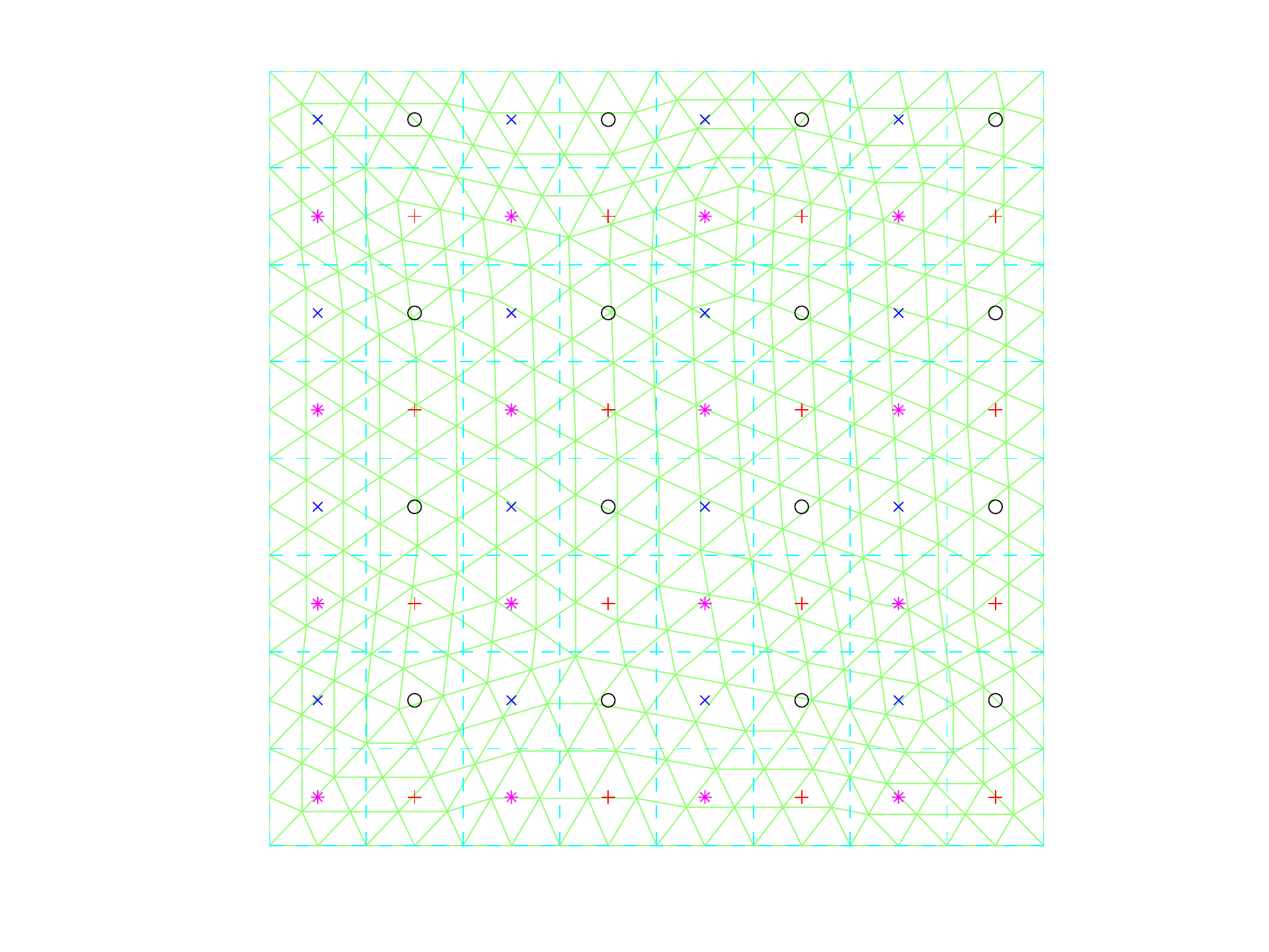}
\caption{Coloring on the finest level $L$}
\label{fig:cgs}
\end{figure}

Thanks to the auxiliary structured grid again, we have only one unstructured grid and it is the finest level, and all the other coarse levels are structured.  Therefore, the coloring for the coarse levels becomes trivial.  Because we have a 9-point stencil, $4$ colors are sufficient for the coarse-level structured grid and a $4$-color point-wise Gauss-Seidel smoother can be used as the parallel smoother on the coarse levels.  For the unstructured grid on the fine level,  we apply a modified coloring strategy.  In stead of coloring each DoF, we color each aggregate.  Because the aggregates are sitting in the auxiliary structured grid which is formed by all the subregions $\Omega_i^L$ on level $L$, $4$ colors are sufficient for coloring.  Therefore, on the finest level $L$, a $4$-color block-wise Gauss-Seidel smoother can be applied with the aggregates serving as nonoverlapping blocks, see Figure \ref{fig:cgs}.

%

\section{Numerical Results}\label{sec:result}
In this section, we perform several numerical experiments to demonstrate the efficiency of our proposed AMG method for GPU implementation.  We test the parallel algorithm on partial differential equations discretized on different grids in 2D or 3D, in order to show the potential of the proposed AMG method for practical applications. 


\subsection{Test Platform}
Our focal test and comparison platform is a NVIDIA Tesla C2070 together with a Dell computing workstation. Details in regard to the machine are given in Table \ref{tab:test_platform}.

\begin{table}[htdp]
\caption{Test Platform}
\begin{center}
\begin{tabular}{r|l}
\hline \hline
CPU Type & Intel  \\
CPU Clock & 2.4 GHz \\
Host Memory & 16 GB  \\
\hline 
GPU Type & NVIDIA Telsa C2070 \\
GPU Clock &  575MHz \\
Device Memory & 6 GB \\
CUDA Capability &   2.0\\
\hline
Operating System & RedHat  \\
CUDA Driver &   CUDA 4.1\\
Host Complier &   gcc 4.1\\
Device Complier &  nvcc 4.1\\
CUSP &  v0.3.0\\
\hline \hline
\end{tabular}
\end{center}
\label{tab:test_platform}
\end{table}

In our numerical tests, because our aim is to demonstrate the improvement of our new algorithm on GPUs,  we concentrate on comparing our proposed method with the parallel smoothed aggregation AMG method implemented in the CUSP package \cite{Garland2010}.  CUSP is an open source C++ library of generic parallel algorithms for sparse linear algebra and graph computations on CUDA-enabled GPUs. It provides a flexible, high-level interface for manipulating sparse matrices and solving sparse linear systems, and all CUSP's algorithms and implementations have been optimized for GPU by NVIDIA's research group.  To the best of our knowledge, the parallel AMG method implemented in the CUSP package is the state-of-the-art AMG method on GPU. 


%

\subsection{Performance}
\begin{example}[2D Poisson problem] \label{exp:1}
Consider the following Poisson equation in 2D,
\begin{equation*}
 - \Delta u = f,  \quad \text{in} \ \Omega \subset \mathbb{R}^2,
\end{equation*}
with the Dirichlet boundary condition
\begin{equation*}
u = 0, \quad \text{on} \ \partial \Omega.
\end{equation*}
\end{example}
The standard linear finite element method is used to solve Example \ref{exp:1} for a certain triangulation of $\Omega$.

\noindent {\bf 2D uniform grid on $\Omega = (0,1) \times (0,1)$.} On a 2D uniform grid, the standard linear finite element method for the Poisson equation gives a 5-point stencil and the resulting stiffness matrix has a banded data structure.  In this case, our proposed aggregation algorithm coincides with the standard geometric coarsening, which suggests that our proposed AMG method has clear advantages over other AMG methods. Table \ref{table: uniform} shows the numbers of iterations and wall time required to reduce the relative residual to $10^{-6}$ for different implementations.

\begin{table}[!htp]
\begin{center}
\caption{Wall time and number of iterations for the Poisson problem on a 2D uniform grid}
\begin{tabular}{|c|c|c|c|c|c|c|c|c|}
\hline
& \multicolumn{4}{|c|}{1024$\times$1024}&\multicolumn{4}{|c|}{2048$\times$2048}\\  \hline
& Setup & \# Iter & Solve &Total &Setup &\# Iter& Solve &Total\\ \hline
CUSP&0.63&36&0.35&0.98&2.38&41&1.60&3.98\\ \hline
New &0.03&10&0.13&0.16&0.11&11&0.43&0.54\\ \hline
\end{tabular}
\label{table: uniform}
\end{center}
\end{table}

As shown in Table \ref{table: uniform}, compared with the CUSP package, our proposed AMG method is 
about $21$ to $22$ times faster in the setup phase, $3$ to $4$ times faster in the solver phase, and $6$ to $7$ time faster in total.  

\begin{figure}[!htp]
\small
\centering
\includegraphics[width=0.7\textwidth]{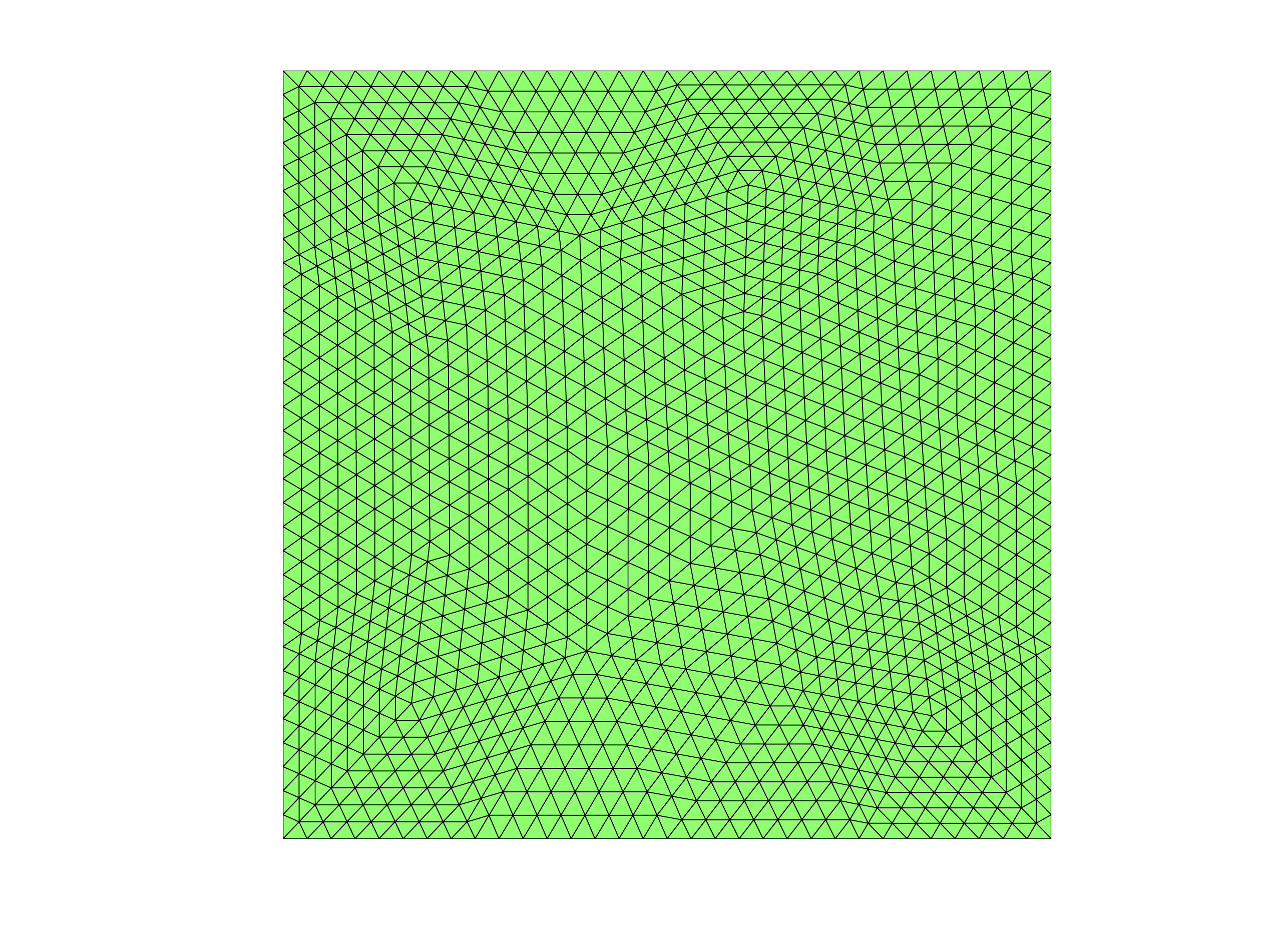}
\caption{Quasi-uniform grid for a 2D unit square}
\label{fig:quasi-uniform}
\end{figure}

\noindent {\bf 2D quasi-uniform grid on $\Omega = (0,1) \times (0,1)$.} The quasi-uniform grid on the unit square is shown in Figure \ref{fig:quasi-uniform}.  Table \ref{label: quasi-uniform} shows the comparison between the CUSP package and our method.  For the problem with $1.4$ million unknowns,  our new method runs $3.03$ times faster than the AMG method in the CUSP package, and for the problem with $5.7$ million unknowns, our method runs $3.82$ times faster than the CUSP package.  The speedup is even more significant in the setup phase: our method is about $7$ times faster than the setup phase implemented in the CUSP package, thus demonstrating the efficiency of our proposed auxiliary-grid-based aggregation algorithm. 

\begin{table}[!htp]
\begin{center}
\caption{Wall time and number of iterations for the Poisson problem on a 2D quasi-uniform grid}
\label{label: quasi-uniform}
\begin{tabular}{|c|c|c|c|c|c|c|c|c|}
\hline
& \multicolumn{4}{|c|}{\# Unknowns = $1,439,745$}&\multicolumn{4}{|c|}{\# Unknowns = $5,763,073$}\\  \hline
& Setup & \# Iter &Solve&Total&Setup& \# Iter &Solve&Total\\ \hline
CUSP&1.09&40&0.82&1.91&3.65&50&4.11&7.76\\ \hline
New&0.16&15&0.49&0.65&0.53&13&1.52&2.05\\ \hline
\end{tabular}
\end{center}
\end{table}

\begin{figure}[!htp]
\small
\centering
\includegraphics[width=0.7\textwidth]{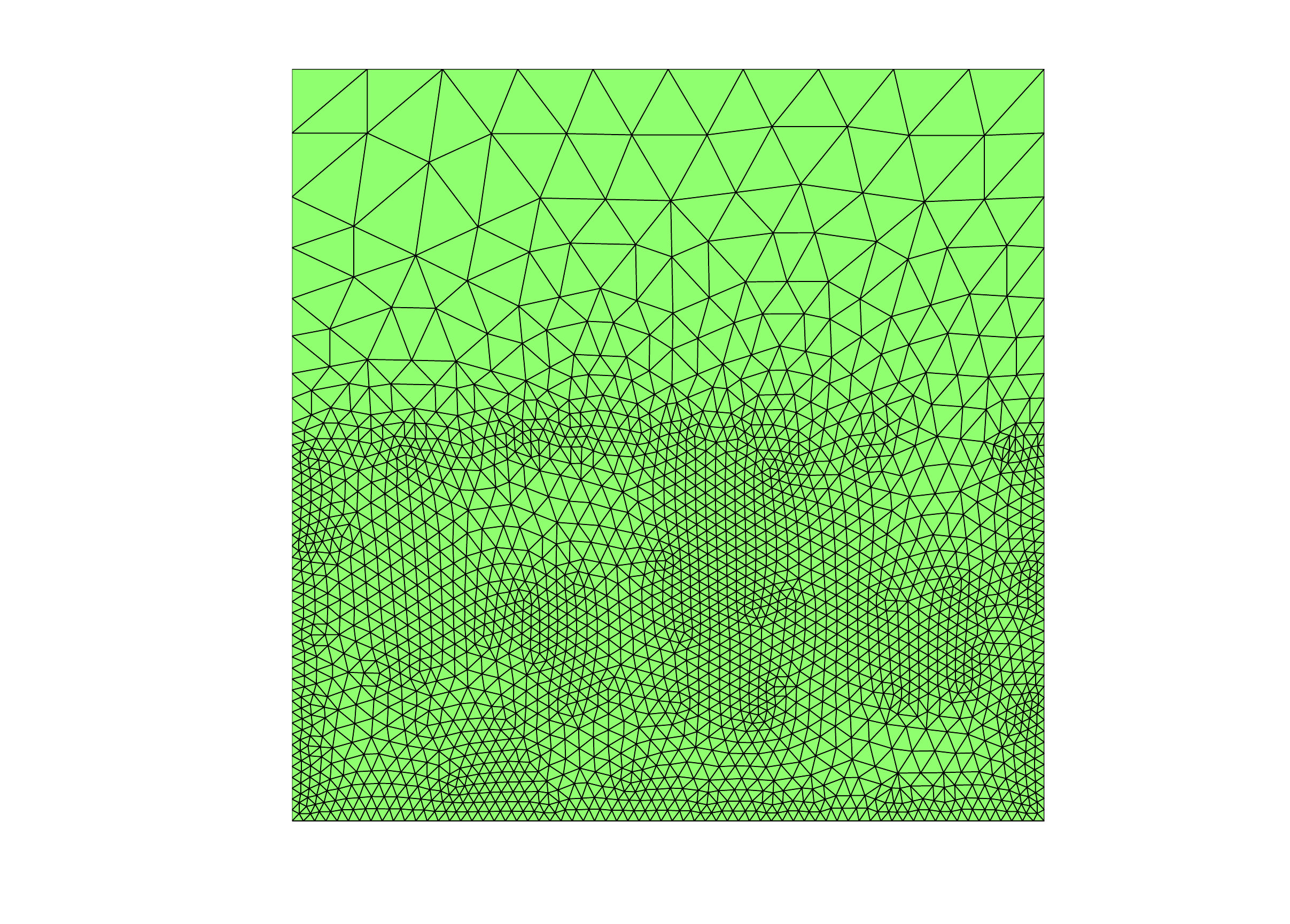}
\caption{Shape-regular grid for a 2D unit square}
\label{fig:shaperegular}
\end{figure}

\noindent {\bf 2D shape-regular grid on $\Omega = (0,1) \times (0,1)$.} We also test our auxiliary-grid-based AMG method for the shape-regular grid (see Figure \ref{fig:shaperegular}).  Table \ref{label:shaperegular} shows the number of iterations and the wall time.  We can see that, compared with the CUSP package, our proposed auxiliary-grid-based AMG method can achieve about $6$ times speedup in the setup phase and $2$ times speedup in total.  The results are not as good as those for the uniform or quasi-uniform grid because the local density of the grid varies and thereby causes an unbalanced distribution of unknowns at the finest level.  However, our method still achieves reasonable speedup in the setup phase, which is considered the most challenging aspect for a parallel AMG method.  However, for a large-scale problem i.e., one with $13$ million unknowns, the CUSP package terminates due to the limitation of the memory in the GPU.  Our method requires much less memory than the CUSP package does and still produces reasonable results in an optimal way. 

\begin{table}[!htp]
\begin{center}
\caption{Wall time and number of iterations for the Poisson problem on a 2D shape-regular grid}
\label{label:shaperegular}
\begin{tabular}{|c|c|c|c|c|c|c|c|c|}
\hline
& \multicolumn{4}{|c|}{\# Unknowns = $3,404,545$}&\multicolumn{4}{|c|}{\# Unknowns = $13,624,833$}\\  \hline
& Setup & \# Iter &Solve&Total&Setup& \# Iter &Solve&Total\\ \hline
CUSP&2.23&52&2.51&4.74& -- &--&--&--\\ \hline
New &0.35&25&2.12&2.47&1.62&27&8.67&10.29\\ \hline
\end{tabular}
\end{center}
\end{table}

\begin{figure}[!htp]
\small
\centering
\includegraphics[width=0.7\textwidth]{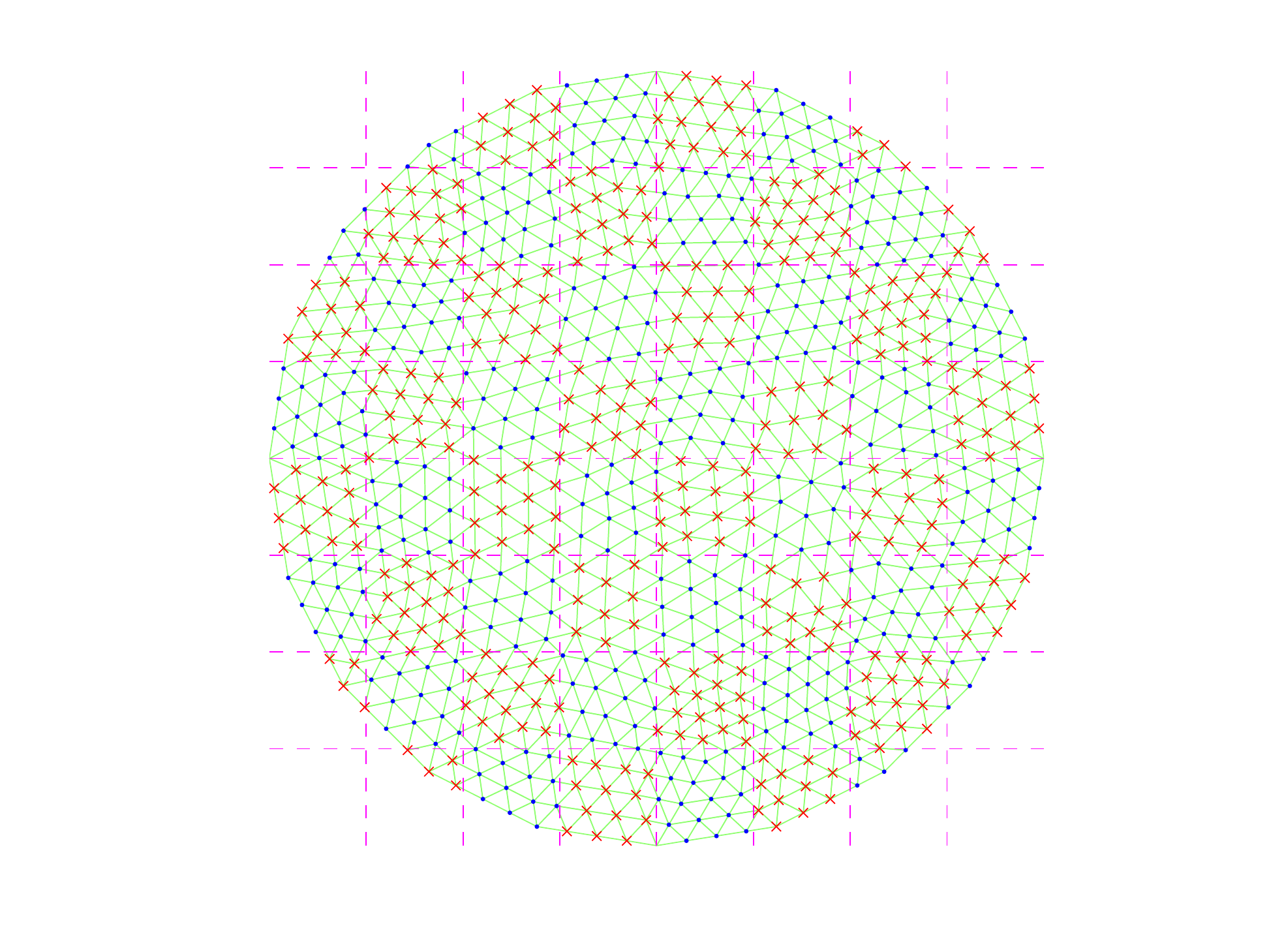}
\caption{Shape-regular grid for a 2D circle domain}
\label{fig:circle}
\end{figure}

\noindent {\bf 2D quasi-uniform grid on a disk domain.} Instead of the square unit domain, we test the performance of our method on a disk domain as shown in Figure \ref{fig:circle}.  Although there might be some empty aggregates that they may affect the overall efficiency, our algorithm is still robust as shown in Table \ref{table:circlegrid}, and can achieve $7$ times speedup in the setup phase and $2$ times speedup in total compared with the AMG method in CUSP.

\begin{table}[!htp]
\begin{center}
\caption{Wall time and number of iterations for the Poisson problem on a disk domain}
\label{table:circlegrid}
\begin{tabular}{|c|c|c|c|c|c|c|c|c|}
\hline
 & \multicolumn{4}{|c|}{\# Unknowns = $2,881,025$}&\multicolumn{4}{|c|}{\# Unknowns = $11,529,733$}\\  \hline
& Setup &\# Iter &Solve&Total&Setup&\# Iter &Solve&Total\\ \hline
CUSP&2.19&50&1.68&3.87&--&--&--&--\\ \hline
New &0.30&21&1.65&1.85&1.37&23&7.21&8.58\\ \hline
\end{tabular}
\end{center}
\end{table}

\begin{example}[3D heat-transfer problem] \label{exp:2}
Consider 
\begin{equation*}
\left\{
\begin{aligned}
\frac{\partial T}{\partial t} - \nabla \cdot (D \nabla T) &= 0, \quad \text{ in } \Omega \in \mathbb{R}^3\\
T &= 356,\quad \text{ on }\Gamma_{inlet}\\
\frac{\partial T}{\partial n} &= 0 \quad \text{ on }  \Gamma_{D}
\end{aligned}
\right.
\end{equation*}
with finite volume discretization, where $T_0=256$, $D = 1$, and the boundary $\Gamma = \Gamma_{inlet}\cup\Gamma_{D}$.
\end{example}

\begin{figure}[!htp]
\small
\centering
\includegraphics[width=0.4\textwidth]{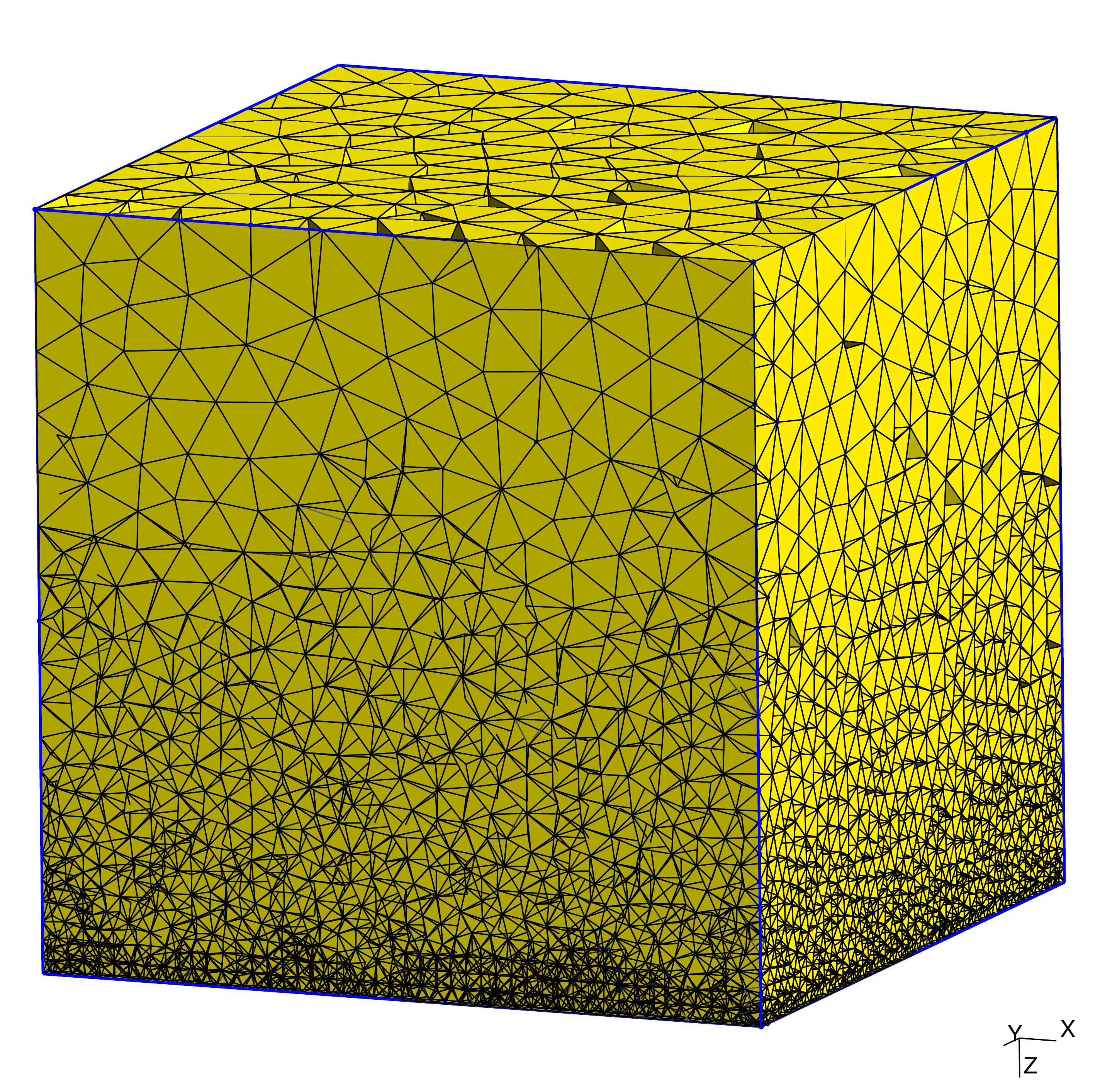}
\caption{Shape regular-grid for the 3D heat transfer problem on a cubic domain}
\label{fig:shaperegular-cube}
\end{figure}

\begin{figure}[!htp]
\small
\centering
\includegraphics[width=0.4\textwidth]{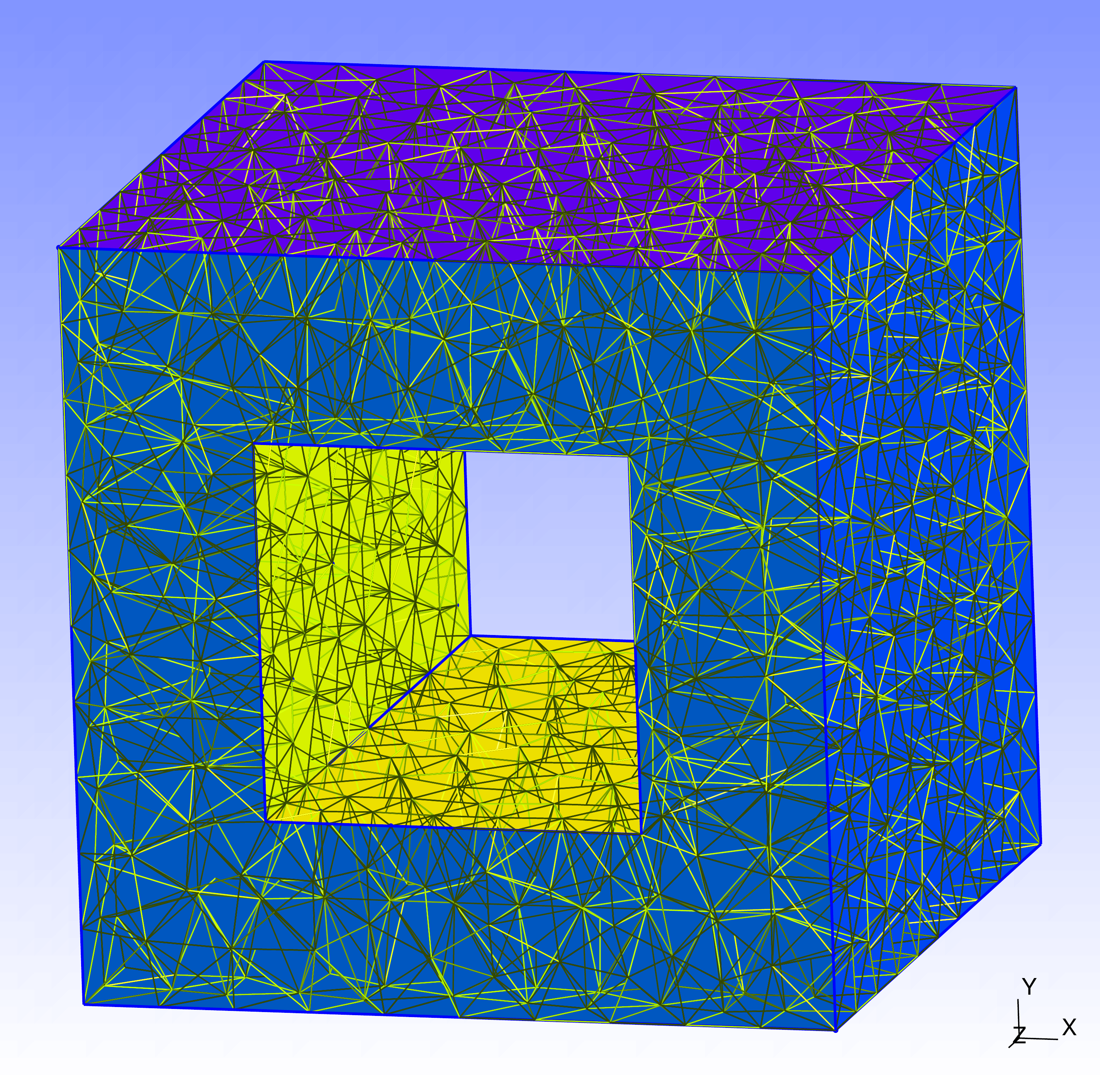}
\caption{Shape regular-grid for the 3D heat transfer problem on a cavity domain}
\label{fig:shaperegular-cs}
\end{figure}

We consider two 3D computational domains, of which one is a unit cube (Figure \ref{fig:shaperegular-cube}) and the other is a cavity domain (Figure \ref{fig:shaperegular-cs}).  Shape-regular grids are used on both domains.  The numerical results for one time step are shown in Tables \ref{label:cube} and \ref{label:concentric square} for the cubic and cavity domain, respectively.  Compared with the other aggregation-based AMG method on GPU, we can see about $6$ times speedup in the setup phase, $2$ times speedup in the solver phase, and about $2$ to $5$ times speedup in total.  The numerical results demonstrate the efficiency and robustness of our proposed AMG algorithm for isotropic problems in 3D. And these results affirm the potential of the new method for solving practical problems.

\begin{table}[!htp]
\begin{center}
\caption{Time/Number of iterations for the heat-transfer problem on a 3D unit cube}
\label{label:cube}
\begin{tabular}{|c|c|c|c|c|c|c|}
\hline
& \multicolumn{3}{|c|}{ \# Unknowns = $909,749$}&\multicolumn{3}{|c|}{\# Unknowns = $1,766,015$}\\  \hline
& Setup & Solve&Total&Setup&Solve&Total\\ \hline
CUSP&1.75&0.60&2.35&1.94&0.52&2.46\\ \hline
New &0.26&0.37&0.63&0.39&0.64&1.03\\ \hline
\end{tabular}
\end{center}
\end{table}

\begin{table}[!htp]
\begin{center}
\caption{Wall time and the number of iterations for the heat-transfer problem on a cavity domain}
\label{label:concentric square}
\begin{tabular}{|c|c|c|c|c|c|c|}
\hline
 & \multicolumn{3}{|c|}{\# Unknowns = $410,482$}&\multicolumn{3}{|c|}{\# Unknowns = $2,094,240$}\\  \hline
& Setup & Solve&Total&Setup&Solve&Total\\ \hline
CUSP&0.67&0.12&0.79&1.14&0.38&1.52\\ \hline
New &0.11&0.05&0.16&0.22&0.32&0.54\\ \hline
\end{tabular}
\end{center}
\end{table}

\section{Conclusion}\label{sec:conclusion}
In this paper, we develop a new parallel auxiliary-grid-based AMG method.  The proposed AMG algorithm is based on an unsmoothed aggregation AMG method and a nonlinear AMLI-cycle (K-cycle) method.  The coarsening and smoothing procedures in our new algorithm are based on a quadtree in 2D (an octree in 3D) generated from an auxiliary structured grid. This provides a (nearly) optimal load balance and predictable communication patterns--factors that make our new algorithm suitable for parallel computing, especially on GPU.  More precisely, the special structure of the auxiliary grid narrows the bandwidth of the coarse-grid matrix (9-point stencil in 2D and 27-point stencil in 3D), which gives us explicitly control of the sparsity pattern of the coarse-grid matrices and reduces the operator complexity of the AMG methods.  Moreover, due to the regular sparsity pattern,  we can use the ELLPACK format, which enables efficient sparse matrix-vector multiplication on GPU. In addition, the auxiliary grid allows us to use a colored Gauss-Seidel smoother without any extra work, thus improving not only the convergence rate but also the parallel performance of the UA-AMG solver.  Numerical results show that our new method can speed up the existing state-of-the-art GPU code (CUSP from NVIDIA) by about $4$ times on a quasi-uniform grid and by $2$ times on a shape-regular grid for certain model problems in both 2D or 3D.


\bibliographystyle{plain}
\bibliography{UAAMGonGPU}

\end{document}